\documentclass[10pt]{article}
\newtheorem{theorem}{Theorem}
\newtheorem{proposition}{Proposition}
\newtheorem{lemma}{Lemma}
\newtheorem{corollary}{Corollary}

\input tcilatex

\begin{document}

\begin{titlepage}

\begin{center}

{\LARGE \bf On the existence of a new type of periodic and quasi-periodic 
orbits for twist maps of the torus }

\end{center}

\vskip  0.5truecm

\centerline {{\large Salvador Addas Zanata}}

\vskip 0.3truecm

\centerline { {\sl Department of Mathematics}}
\centerline {{\sl Princeton University}}
\centerline {{\sl Fine Hall-Washington Road, Princeton, NJ 08544-1000,
USA}}

\vskip .3truecm

\begin{abstract}
We prove that for a large and important class of $C^1$
twist maps of the torus periodic and quasi-periodic orbits of a new type
 exist, provided that there are no rotational invariant circles
(R.I.C's). These orbits have a non-zero ''vertical rotation number'' (V.R.N.), 
in contrast to what happens to Birkhoff periodic orbits
and Aubry-Mather sets. The V.R.N. is rational for a periodic orbit and
irrational for a quasi-periodic. We also prove that the existence of
an orbit with a $V.R.N=a>0,$ implies the existence of orbits with $V.R.N=b,$
for all $0<b<a.$ In this way, related to a generalized definition of rotation number,
 we characterize all kinds of 
periodic and quasi-periodic orbits a twist map of the torus can have.
And as a consequence of the previous
results we obtain that a twist map of the torus with no R.I.C's has positive topological
 entropy, which is a very classical result. 
In the end of the paper we present some examples, 
like the Standard map, such that our results apply.
\end{abstract} 

\vskip .3truecm


\hskip1.8truecm vertical rotation number, Nielsen-Thurston theory

\vskip 0.5truecm
\noindent {\bf E-mail}: szanata@math.princeton.edu

\vfill
\hrule
\noindent{\footnotesize{supported by CNPq, grant number: 200564/00-5 (part
of this work was done while the author was under support by FAPESP,
grant number: 96/08981-3)}}

\end{titlepage}

\baselineskip=6.5mm

\section{Introduction and statements of the principal results}

Twist maps are $C^1$ diffeomorphisms of the cylinder (or annulus, torus)
onto itself that have the following property: The angular component of the
image of a point increases as the radial component of the point increases
(more precise definitions will be given below). Such maps were first studied
in connection with the three body problem by Poincar\'e and later it was
found that first return maps for many problems in Hamiltonian dynamics are
actually twist maps. Although they have been extensively studied, there are
still many open questions about their dynamics. A great progress has been
achieved in the nearly integrable case, by means of KAM theory (see \cite
{moser}) and many important results have been proved in the general case,
concerning the existence of periodic and quasi-periodic orbits\
(Aubry-Mather sets), see \cite{mather}, \cite{aubry}, \cite{katok}. In this
work a result is proved associating the non-existence of rotational
invariant circles to the appearance of periodic and quasi-periodic orbits of
a new type for an important class of twist maps of the torus.

\begin{description}
\item[Notation and definitions]  :

0) Let $(\phi ,I)$ denote the coordinates for the cylinder $S^1\times ${\rm I%
}\negthinspace {\rm R=(}${\rm I}\negthinspace {\rm R/Z\negthinspace
\negthinspace Z}${\rm )}$\times {\rm I}\negthinspace {\rm R}${\rm ,} where $%
\phi $ is defined modulo 1{\rm .} Let $(\widetilde{\phi },\widetilde{I})$
denote the coordinates for the universal cover of the cylinder, {\rm I}%
\negthinspace {\rm R}$^{{\rm 2}}.$ 
For all maps $\widehat{f}:S^1\times ${\rm I}\negthinspace {\rm R}$%
\rightarrow S^1\times ${\rm I}\negthinspace {\rm R} we define 
$$
\begin{array}{c}
(\phi ^{\prime },I^{\prime })=
\widehat{f}(\phi ,I) \\ and \\ 
(\widetilde{\phi }^{\prime },\widetilde{I}^{\prime })=f(\widetilde{\phi },
\widetilde{I})
\end{array}
, 
$$
where $f:{\rm I}\negthinspace {\rm R^2\rightarrow I}\negthinspace {\rm R^2}$
is a lift of $\widehat{f}.$

1) $D_r^1({\rm I}\negthinspace {\rm R^2)}=\{f:{\rm I}\negthinspace {\rm %
R^2\rightarrow I}\negthinspace {\rm R^2}$ / $f$ is a $C^1$ diffeomorphism of
the plane, $\widetilde{I}^{\prime }(\widetilde{\phi },\widetilde{I})%
\stackrel{\widetilde{I}\rightarrow \pm \infty }{\rightarrow }\pm \infty $, $%
\partial _{\widetilde{I}}\widetilde{\phi }^{\prime }$$>0$ (twist to the
right), $\widetilde{\phi }^{\prime }(\widetilde{\phi },\widetilde{I})%
\stackrel{\widetilde{I}\rightarrow \pm \infty }{\rightarrow }\pm \infty $
and $f$ is the lift of a $C^1$ diffeomorphism $\widehat{f}:S^1\times ${\rm I}%
\negthinspace {\rm R}$\rightarrow S^1\times ${\rm I}\negthinspace 
{\rm R}${\rm \}}$.

2) $Diff_r^1(S^1\times {\rm I}\negthinspace {\rm R})$=$\{$ $\widehat{f}%
:S^1\times ${\rm I}\negthinspace {\rm R}$\rightarrow S^1\times ${\rm I}%
\negthinspace 
{\rm R }/ $\widehat{f}$ is induced by an element of $D_r^1$({\rm I}%
\negthinspace {\rm R$^2$)}$\}.$

3) Let $p_1:{\rm I}\negthinspace {\rm R^2\rightarrow I}\negthinspace {\rm R}$
and $p_2:{\rm I}\negthinspace {\rm R^2\rightarrow I}\negthinspace {\rm R}$
be the standard projections, respectively in the $\widetilde{\phi }$ and $
\widetilde{I}$ coordinates ($p_1(\widetilde{\phi },\widetilde{I})=\widetilde{%
\phi }$ and $p_2(\widetilde{\phi },\widetilde{I})=\widetilde{I}$). We also
use $p_1$ and $p_2$ for the standard projections of the cylinder.

4) Given a measure $\mu $ on the cylinder that is positive on open sets,
absolutely continuous with respect to the Lebesgue measure and a map $
\widehat{T}\in Diff_r^1(S^1\times {\rm I}\negthinspace {\rm R})$ we say that 
$\widehat{T}$ is $\mu $-exact if $\mu $ is invariant under $\widehat{T}$ and
for any open set $A$ homeomorphic to the cylinder we have: 
\begin{equation}
\label{miexact}\mu (\widehat{T}(A)\backslash A)=\mu (A\backslash \widehat{T}%
(A))
\end{equation}
For an area-preserving map $\widehat{T}\in Diff_r^1(S^1\times {\rm I}%
\negthinspace {\rm R}),$ there is a simple criteria to know if it is exact
or not. $\widehat{T}$ is exact if and only if its generating function $h(
\widetilde{\phi },\widetilde{\phi }^{\prime }),$ defined on ${\rm I}%
\negthinspace {\rm R^2,}$ satisfies $h(\widetilde{\phi }+1,\widetilde{\phi }%
^{\prime }+1)=h(\widetilde{\phi },\widetilde{\phi }^{\prime })\ ($see$\ \cite
{meiss}).$


5) Let $TQ\subset D_r^1({\rm I}\negthinspace {\rm R^2})$ be such that for
all $T\in TQ$ we have:

\begin{equation}
\label{DefTQ} 
\begin{array}{c}
T:\left\{ 
\begin{array}{c}
\widetilde{\phi }^{\prime }=T_\phi (\widetilde{\phi },\widetilde{I}) \\ 
\widetilde{I}^{\prime }=T_I(\widetilde{\phi },\widetilde{I}) 
\end{array}
\right. , 
\text{ with }\partial _{\widetilde{I}}\widetilde{\phi }^{\prime }=\partial
_{ \widetilde{I}}T_\phi (\widetilde{\phi },\widetilde{I})>0\text{ and:} \\ 
\left\{ 
\begin{array}{l}
T_I( 
\widetilde{\phi }+1,\widetilde{I})=T_I(\widetilde{\phi },\widetilde{I}) \\ 
T_I( 
\widetilde{\phi },\widetilde{I}+1)=T_I(\widetilde{\phi },\widetilde{I})+1 \\ 
\text{ }T_\phi (\widetilde{\phi }+1,\widetilde{I})=T_\phi (\widetilde{\phi }%
, \widetilde{I})+1\text{ } \\ T_\phi (\widetilde{\phi },\widetilde{I}%
+1)=T_\phi (\widetilde{\phi },\widetilde{I})+1 
\end{array}
\right. 
\end{array}
\end{equation}

Every $T\in TQ$ induces a map $\widehat{T}\in Diff_r^1(S^1\times {\rm I}%
\negthinspace {\rm R})$ and a map $\overline{T}:{\rm T^2\rightarrow T^2,}$
where ${\rm T^2=I\negthinspace R^2/Z\negthinspace \negthinspace Z^2}$ is the
2-torus. Let $p:{\rm I\negthinspace R^2\rightarrow T^2}$ be the associated
covering map.

6) Given $T\in TQ,$ we say that $\beta \in ]0,\pi /2[$ is a uniform angle of
deviation for $T$ if 
$$
DT\mid _x.\left( 
\begin{array}{c}
0 \\ 
1 
\end{array}
\right) \in C_I(\beta )\text{ and }DT^{-1}\mid _z.\left( 
\begin{array}{c}
0 \\ 
1 
\end{array}
\right) \in C_{II}(\beta ), 
$$
for all $x,z\in {\rm I\negthinspace R^2,}$ where $C_I(\beta )$ and $%
C_{II}(\beta )$ are the angular sectors: 
$$
\begin{array}{c}
C_I(\beta )=\left\{ (\phi ,I)\in 
{\rm I\negthinspace R^2:}\ \phi >0\ \text{and}\ -\text{cotan}(\beta ).\phi
\leq I\leq \text{cotan}(\beta ).\phi \right\} \\ C_{II}(\beta )=\left\{
(\phi ,I)\in {\rm I\negthinspace R^2:}\ \phi <0\ \text{and}\ \text{cotan}%
(\beta ).\phi \leq I\leq -\text{cotan}(\beta ).\phi \right\} 
\end{array}
$$

7) Let $\pi :{\rm I\negthinspace R^2\rightarrow S^1\times I\negthinspace R}$
be the following covering map:

\begin{equation}
\label{piK}\pi (\widetilde{\phi },\widetilde{I})=\left( \widetilde{\phi }\
(mod\ 1),\ \widetilde{I}\right) 
\end{equation}

8) For maps of the torus we can generalize the notion of rotation number,
originally defined for circle homeomorphisms, as follows:

Given a map $\overline{f}:{\rm T^2\rightarrow T^2}$ and $x\in {\rm T^2,}$
let $f:{\rm I\negthinspace R^2\rightarrow I\negthinspace R^2}$ be a lift of $
\overline{f}$ and $\widetilde{x}\in p^{-1}(x).$ The rotation vector $\rho
(x,f)$ is defined as (let $\overline{{\rm I\negthinspace R}}={\rm I%
\negthinspace R\cup \{-\infty \}\cup \{+\infty \}}$)

\begin{equation}
\label{rotnum}\rho (x,f)=\stackunder{n\rightarrow \infty }{\lim }\frac{f^n( 
\widetilde{x})-\widetilde{x}}n\in \overline{{\rm I\negthinspace R}}{\rm ^2,} 
\text{ if the limit exists.} 
\end{equation}
\end{description}

Of course for a map $\overline{f}:{\rm T^2\rightarrow T^2}$ that is not
homotopic to the identity${\rm ,}$ the above limit may depend on the choice
of $\widetilde{x}\in p^{-1}(x).$ So, for all $T\in TQ$ we have to modify a
little the above definition for rotation vector.\ The following lemma shows
what type of changes must be done:

\begin{lemma}
\label{characrot}: Given $T\in TQ$ and $\widetilde{z}\in {\rm I\negthinspace %
R^2}$ such that $\rho \left( \widetilde{z},T\right) =\stackunder{%
n\rightarrow \infty }{\lim }\dfrac{T^n(\widetilde{z})-\widetilde{z}}n$
exists, we have two possibilities:

1) $\rho (\widetilde{z},T)=(\omega ,0),$ with $\omega \in {\rm I%
\negthinspace R,}$

2) $\rho (\widetilde{z},T)=(+\infty ,\omega )$ or $(-\infty ,\omega ),$ with 
$\omega \in {\rm I\negthinspace R.}$
\end{lemma}

Given $i,j\in {\rm Z\negthinspace \negthinspace Z,}$ we have:

$\bullet $ $\rho (\widetilde{z},T)=(\omega ,0)$ $\Rightarrow $ $\rho \left( 
\widetilde{z}+(i,j),T\right) =(\omega +j,0)$

$\bullet $ $\rho (\widetilde{z},T)=(\pm \infty ,\omega )$ $\Rightarrow $ $%
\rho \left( \widetilde{z}+(i,j),T\right) =(\pm \infty ,\omega )$

So, for all $\widetilde{x}\in {\rm I\negthinspace R^2}$ such that $\rho ( 
\widetilde{x},T)$ exists, there are 2 different cases:

$\bullet $ Case 1: $p_1\circ \rho (\widetilde{x},T)\in {\rm I\negthinspace R.%
}$ In this case we shall define the rotation vector of $x=p(\widetilde{x}%
)\in {\rm T^2}$ as follows:

\begin{equation}
\label{rotnum1} \rho (x,T)=\left( \stackunder{n\rightarrow \infty }{\lim } 
\dfrac{p_1\circ T^n(\widetilde{x})-p_1(\widetilde{x})}n\ (mod\ 1),\text{ }%
0\right) 
\end{equation}

$\bullet $ Case 2: $\left| p_1\circ \rho (\widetilde{x},T)\right| =\infty .$
In this case $\rho (x,T)$ is defined as in (\ref{rotnum}): 
\begin{equation}
\label{rotnum2}\rho (x,T)=\stackunder{n\rightarrow \infty }{\lim }\dfrac{%
T^n( \widetilde{x})-\widetilde{x}}n,\text{ where }x=p(\widetilde{x})\text{ } 
\end{equation}

We just remark that even in this case $p_2\circ \rho (x,T)$ may be zero.

When $f:{\rm I\negthinspace R^2\rightarrow I\negthinspace R^2}$ induces a
map $\overline{f}:{\rm T^2\rightarrow T^2}$ homotopic to the identity map
(in this case the rotation vector is given by expression (\ref{rotnum})), in
many situations (see \cite{franks}, \cite{LMack} and \cite{misiu}) we can
guarantee the existence of a convex open set $B\subset {\rm I\negthinspace %
R^2,}$ such that for all $v\in B,$ $\exists $ $x\in {\rm T^2}$ such that $%
\rho (x,f)=v$ and if $v=(\frac rq,\frac sq),$ then $x\in {\rm T^2}$ can be
chosen such that $f^q(x)=x+(r,s).$ A major difference in the situation
studied here is that given $T\in TQ,$ as we have already said, the
diffeomorphism $\overline{T}:{\rm T^2\rightarrow T^2}$ induced by $T$ is not
homotopic to the identity. In fact, it is homotopic to the following linear
map (where $\phi $ and $I$ are taken $mod$ $1$):

\begin{equation}
\label{represent} \left( 
\begin{array}{c}
\phi ^{\prime } \\ 
I^{\prime } 
\end{array}
\right) =\left( 
\begin{array}{cc}
1 & 1 \\ 
0 & 1 
\end{array}
\right) .\left( 
\begin{array}{c}
\phi \\ 
I 
\end{array}
\right) 
\end{equation}

Before presenting the first theorem we still need more definitions.

\begin{description}
\item[Definitions]  : Given $T\in TQ,$ let $\overline{T}:{\rm T^2\rightarrow
T^2}$ be the torus diffeomorphism induced by $T.$

$\bullet $ We say that $x\in {\rm T^2}$ belongs to a $n-$periodic orbit (or
set), if for some $n\in {\rm I\negthinspace N^{*}}$ we have $\overline{T}%
^n(x)=x$ and for all $m\in {\rm I\negthinspace N^{*},}$ $0<m<n,$ $\overline{T%
}^m(x)\neq x.$ So the periodic orbit which $x$ belongs is $O_x=\{x,\overline{%
T}(x),...,\overline{T}^{n-1}(x)\}.$ In this case we have the following
implications (now we consider the standard projections $p_i:\overline{{\rm I}%
\negthinspace {\rm R}}{\rm ^2\rightarrow }\overline{{\rm I}\negthinspace 
{\rm R}},$ $i=1,2$): 
$$
\left\{ 
\begin{array}{l}
p_2\circ \rho (x,T)=0 
\text{ }\Rightarrow \text{ }p_1\circ \rho (x,T)\text{ is a rational number,}
\\ p_1\circ \rho (x,T)=\pm \infty \text{ }\Rightarrow \text{ }p_2\circ \rho
(x,T)\text{ is a non-zero rational number.} 
\end{array}
\right. 
$$
$\bullet $ We say that $Q\subset {\rm T^2}$ is a quasi-periodic set for $
\overline{T}$ in the following cases:

1) $Q$ is the projection of an Aubry-Mather set $\widehat{Q}\subset {\rm %
S^1\times I\negthinspace R.}$ In this case $p_1\circ \rho (z,T)\in [0,1[$ is
an irrational number which does not depend on the choice of $z\in Q.$

2) $Q$ is a compact $\overline{T}-invariant$ set such that $p_2\circ \rho
(z,T)$ is an irrational number which does not depend on the choice of $z\in
Q $.
\end{description}

As before, let $T\in TQ$ and $\overline{T}:{\rm T^2\rightarrow T^2}$ be the
diffeomorphism induced by $T.$ So as a simple consequence of lemma (\ref
{characrot}) we have the following classification theorem:

\begin{theorem}
\label{orbtype}: Let $x\in {\rm T^2}$ belong to a periodic or a
quasi-periodic set. Then there are two different situations:

1) $\exists $ $C>0$ such that $\left| p_2\circ T^n(\widetilde{x})-p_2( 
\widetilde{x})\right| <C,$ for all $n>0$ and $\widetilde{x}\in p^{-1}(x)$ $%
\Rightarrow $ $\rho (x,T)=(\omega ,0)$ for some $\omega \in [0,1[$

2) $p_2\circ T^n(\widetilde{x})\stackrel{n\rightarrow \infty }{\rightarrow }%
\pm \infty ,$ for all $\widetilde{x}\in p^{-1}(x)$ $\Rightarrow $ $\rho
(x,T)=(\pm \infty ,\omega )$ for some $\omega \in {\rm I\negthinspace R^{*}}$%
.
\end{theorem}

Remarks:

$\bullet $ If we call $\widehat{T}:S^1\times {\rm I}\negthinspace {\rm %
R\rightarrow }S^1\times {\rm I}\negthinspace {\rm R}$ the map of the
cylinder induced by $T$, it is easy to see that case 1 above corresponds to
periodic and quasi-periodic orbits for $\widehat{T}.$ These are the standard
periodic and quasi-periodic orbits, whose existence is assured by theorem (%
\ref{aubmath}) (see \cite{mather}, \cite{aubry}, \cite{katok}, \cite
{lecalvez1}, \cite{lecalvez2}).

$\bullet $ Case 2 corresponds to orbits for $\widehat{T}$ that either go up
or down on the cylinder, depending on the sign of $\omega \in {\rm I%
\negthinspace
R^{*}.}$ If $\omega >0$ $(<0)$ then $\rho (x,T)=(+\infty (-\infty ),\omega
). $

$\bullet $ As we said, there may be a point $x\in {\rm T^2}$ such that $\rho
(x,T)=(\pm \infty ,0).$ It is clear that $x$ is not periodic because its $
\overline{T}$-orbit can not be finite and $x$ does not belong to a
quasi-periodic set, because any component of $\rho (x,T)$ is irrational.

A periodic or quasi-periodic orbit $O$ for $\overline{T}$ that belongs to
case 2 in theorem (\ref{orbtype}) can (as the orbits belonging to case 1) be
characterized by a single number, the ''vertical rotation number'', which is
defined in the following way:

\begin{equation}
\label{vertrot}\rho _V(O)=p_2\circ \rho (x,T)=\stackunder{n\rightarrow
\infty }{\lim }\frac{p_2\circ T^n(x)-p_2(x)}n,\text{ for any }x\in O. 
\end{equation}


As we want to characterize all kinds of periodic and quasi-periodic orbits a
twist map can have, we recall a well-known result:

\begin{theorem}
\label{aubmath}: Given a map $T\in TQ$ such that $\widehat{T}:S^1\times {\rm %
I}\negthinspace {\rm R\rightarrow }S^1\times {\rm I}\negthinspace {\rm R}$
is $\mu $-exact for some measure $\mu $ we have: For every $\omega \in {\rm I%
\negthinspace R}$ there is a $\widehat{T}-$periodic or quasi-periodic orbit $%
O\subset S^1\times {\rm I}\negthinspace 
{\rm R,}$ respectively, for rational and irrational values of $\omega ,$
such that $\rho (O)=\omega .$
\end{theorem}

See \cite{mather}, \cite{aubry}, \cite{katok}, \cite{lecalvez1}, \cite
{lecalvez2} for different proofs. And we have the following corollary:

\begin{corollary}
\label{case1}: Given a map $T\in TQ$ such that $\widehat{T}$ is $\mu $-exact
for some measure $\mu $ we have: For every $\omega \in [0,1[$ there is a $
\overline{T}-$periodic or quasi-periodic orbit $O\subset {\rm T^2,}$
respectively, for rational and irrational values of $\omega $, such that $%
\rho (x,T)=(\omega ,0),$ for all $x\in O.$
\end{corollary}

So the first type of orbit that appears in theorem (\ref{orbtype}) always
exists.

Before presenting the next results we need another definition.

\begin{description}
\item[Definition]  : Given a map $T\in TQ$ such that $\widehat{T}$ is $\mu $%
-exact we say that $C$ is a rotational invariant circle (R.I.C.) for $T$ if $%
C$ is a homotopically non-trivial simple closed curve on the cylinder and $
\widehat{T}(C)=C$.
\end{description}

By a theorem essentially due to Birkhoff, $C$ is the graph of some Lipschitz
function $\psi :S^1\rightarrow {\rm I}\negthinspace {\rm R}$ (see \cite
{Katok}, page 430).

The following theorems are the main results of this paper:

\begin{theorem}
\label{PtosPerQuoc}: Let $T\in TQ$ be such that $\widehat{T}$ is $\mu $%
-exact. Then:

- given $k\in {\rm Z\negthinspace \negthinspace Z}^{*},$ $\exists $ $N>0,$
such that $\overline{T}$ has a periodic orbit with $\rho _V$ (vertical
rotation number) $=\frac kN,$ if and only if, $T$ does not have R.I.C's.
\end{theorem}

The next theorem shows how these periodic orbits appear:

\begin{theorem}
\label{ImplicPeriod}: Again, for all $T\in TQ$ such that $\widehat{T}$ is $%
\mu $-exact, if $\overline{T}$ has a periodic orbit with $\rho $$_V=\frac
kN, $ then for every pair $(k^{\prime },N^{\prime })\in {\rm Z%
\negthinspace 
\negthinspace Z}^{*}\times {\rm I}\negthinspace {\rm N}^{*},$ such that $%
0<\left| \frac{k^{\prime }}{N^{\prime }}\right| <\left| \frac kN\right| $
and $k.k^{\prime }>0,$ $\overline{T}$ has at least 2 periodic orbits with
vertical rotation number $\rho _V^{\prime }=\frac{k^{\prime }}{N^{\prime }}.$
\end{theorem}

About the quasi-periodic orbits we have the following:

\begin{theorem}
\label{ImpQper}: For all $T\in TQ$ such that $\widehat{T}$ is $\mu $-exact
we have:

If $\overline{T}$ has an orbit with $\rho _V$$=\omega ,$ then for all $%
\omega ^{\prime }\in {\rm I}\negthinspace {\rm R\backslash }${\footnotesize 
{\sf l}}${\rm \negthinspace
\negthinspace }{\sf Q}$ such that $0<\left| \omega ^{\prime }\right| <\left|
\omega \right| $ and $\omega .\omega ^{\prime }>0$$,$ $\overline{T}$ has a
quasi-periodic set with vertical rotation number $\rho _V^{\prime }=\omega
^{\prime }.$
\end{theorem}

As a consequence of the proof of theorem (\ref{ImpQper}) we prove the
following classical result:

\begin{theorem}
\label{entropia}: Every $T\in TQ$ without R.I.C's such that $\widehat{T}$ is 
$\mu $-exact induces a map $\overline{T}:{\rm T^2\rightarrow T^2}$ such that 
$h(\overline{T})>0,$ where $h(\overline{T})$ is the topological entropy of $
\overline{T}.$
\end{theorem}

Theorem (\ref{orbtype}) is an immediate consequence of lemma (\ref{characrot}%
), which is proved using simple ideas and the structure of the set $TQ.$
Theorems (\ref{PtosPerQuoc}), (\ref{ImplicPeriod}) and (\ref{ImpQper}) are
proved using topological ideas, essentially due to the twist condition and
some results due to Le Calvez (see \cite{lecalvez1}, \cite{lecalvez2} and
the next section). In the proofs of theorems (\ref{ImpQper}) and (\ref
{entropia}) we also use some results from the Nielsen-Thurston theory of
classification of homeomorphisms of surfaces up to isotopy, to isotope the
map to a pseudo-Anosov one and then some results due to M. Handel, to prove
the existence of quasi-periodic orbits with irrational vertical rotation
number.

\section{Basic tools}

First we recall some topological results for twist maps essentially due to
Le Calvez (see \cite{lecalvez1} and \cite{lecalvez2}). Let $\widehat{f}\in
Diff_{r}^{1}(S^{1}\times {\rm I}\!{\rm R})$ and $f\in D_{r}^{1}({\rm I}\!%
{\rm R^{2})}$ be its lifting. For every pair $(p,q),$ $p\in {\rm Z\!\!Z}$
and $q\in {\rm I\!N^{*}}$ we define the following sets:

\begin{equation}
\label{Kpq} 
\begin{array}{c}
\widetilde{K}(p,q)=\left\{ (\widetilde{\phi },\widetilde{I})\in {\rm I}\!%
{\rm R^{2}}\text{: }p_{1}\circ f^{q}(\widetilde{\phi },\widetilde{I})= 
\widetilde{\phi }+p\right\} \\ \text{ and } \\ K(p,q)=\pi \circ \widetilde{K}%
(p,q) 
\end{array}
\end{equation}

Then we have the following:

\begin{lemma}
\label{Lcal1}: For every $p\in {\rm Z\negthinspace \negthinspace Z}$ and $%
q\in {\rm I\negthinspace N^{*},\ }K(p,q)\supset C(p,q),$ a connected compact
set that separates the cylinder.
\end{lemma}

\begin{lemma}
\label{inter}: Let $\widehat{f}\in Diff_r^1(S^1\times {\rm I}\negthinspace 
{\rm R})$ be a $\mu $-exact map. Then the following intersection holds: $
\widehat{f}(C(p,q))\cap C(p,q)\neq \emptyset .$
\end{lemma}

Now we need a few definitions:

For every $q\geq 1$ and $\overline{\phi }\in {\rm I}\negthinspace {\rm R\ }$%
let

\begin{equation}
\label{miq}\mu _{q}(t)=f^{q}(\overline{\phi },t),\text{ for }t\in {\rm I}\!%
{\rm R} 
\end{equation}

We say that the first encounter between $\mu _q$ and the vertical line
through some $\phi _0\in {\rm I}\negthinspace {\rm R}$ is for: 
$$
\begin{array}{c}
t_F\in 
{\rm I}\negthinspace {\rm R}\text{ such that} \\ t_F=\min \{t\in {\rm I}%
\negthinspace {\rm R}\text{: }p_1\circ \mu _q(t)=\phi _0\} 
\end{array}
$$
And the last encounter is defined in the same way: 
$$
\begin{array}{c}
t_L\in 
{\rm I}\negthinspace {\rm R}\text{ such that} \\ t_L=\max \{t\in {\rm I}%
\negthinspace {\rm R}\text{: }p_1\circ \mu _q(t)=\phi _0\} 
\end{array}
$$
Of course we have $t_F\leq t_L.$

\begin{lemma}
\label{Lcal2}: For all $\phi _0,\overline{\phi }\in {\rm I}\negthinspace 
{\rm R,}$ let $\mu _q(t)=f^q(\overline{\phi },t),$ as in (\ref{miq}). So we
have the following inequalities: $p_2\circ \mu _q(t_L)\leq p_2\circ \mu _q( 
\overline{t})\leq p_2\circ \mu _q(t_F),$ for all $\overline{t}\in {\rm I%
\negthinspace R}$ such that $p_1\circ \mu _q(\overline{t})=\phi _0.$
\end{lemma}

For all $s\in {\rm Z\!\!Z}$ and $N\in {\rm I\!N^{*}}$ we can define the
following functions\ on $S^{1}$: 
$$
\begin{array}{c}
\mu ^{-}(\phi )=\min \{p_{2}(Q) 
\text{: }Q\in K(s,N)\text{ and }p_{1}(Q)=\phi \} \\ \mu ^{+}(\phi )=\max
\{p_{2}(Q)\text{: }Q\in K(s,N)\text{ and }p_{1}(Q)=\phi \} 
\end{array}
$$
And we can define similar functions for $\widehat{f}^{N}(K(s,N))$: 
$$
\begin{array}{c}
\nu ^{-}(\phi )=\min \{p_{2}(Q) 
\text{: }Q\in \widehat{f}^{N}\circ K(s,N)\text{ and }p_{1}(Q)=\phi \} \\ \nu
^{+}(\phi )=\max \{p_{2}(Q)\text{: }Q\in \widehat{f}^{N}\circ K(s,N) \text{
and }p_{1}(Q)=\phi \} 
\end{array}
$$

\begin{lemma}
\label{graphu}: Defining Graph\{$\mu ^{\pm }$\}=\{$(\phi ,\mu ^{\pm }(\phi
)):\phi \in S^1$\} we have: 
$$
Graph\{\mu ^{-}\}\cup Graph\{\mu ^{+}\}\subset C(s,N) 
$$

So for all $\phi \in S^1$ we have $(\phi ,\mu ^{\pm }(\phi ))\in C(s,N).$
\end{lemma}

And we have the following simple corollary to lemma (\ref{Lcal2}):

\begin{corollary}
\label{ofpre}: $\widehat{f}^N(\phi ,\mu ^{-}(\phi ))=(\phi ,\nu ^{+}(\phi ))$
and $\widehat{f}^N(\phi ,\mu ^{+}(\phi ))=(\phi ,\nu ^{-}(\phi )).$
\end{corollary}

Now we are going to present a lemma due to M. Casdagli (see \cite{casdagli}%
), that together with lemma (\ref{inter}) guarantees the existence of
periodic orbits with all rational rotation numbers, for all $\mu $-exact $
\widehat{f}\in Diff_r^1(S^1\times {\rm I}\negthinspace {\rm R})$.

\begin{lemma}
\label{allrotnum}: If $z\in C(s,N)\cap \widehat{f}(C(s,N))\Rightarrow z$ is $%
(s,N)$ periodic for $\widehat{f}.$
\end{lemma}

We say that $z$ is $(s,N)$ periodic for $\widehat{f}$ if 
$$
\widehat{f}^N(z)=z\ \text{and }\frac{p_1\circ f^N(\widetilde{z})-p_1( 
\widetilde{z})}N=\frac sN\text{ ,} 
$$
where $f:{\rm I}\negthinspace {\rm R^2\rightarrow I}\negthinspace {\rm R^2}$
is a lift of $\widehat{f}$ and $\widetilde{z}\in \pi ^{-1}(z).$

For proofs of all the previous results see Le Calvez \cite{lecalvez1} and 
\cite{lecalvez2}.

The following is another classical result (due to Birkhoff) with some small
changes:

\begin{theorem}
\label{orbclimb}: Given $T\in TQ$ without R.I.C's such that $\widehat{T}$ is 
$\mu $-exact we have:

For all $s,l\in {\rm Z\negthinspace \negthinspace Z,}$ $s>0$ and $l<0,$ $%
\exists $ $P,Q\in S^1\times [0,1]$ and numbers

$1<n_P,n_Q\in {\rm I\negthinspace N}$ such that $\left\{ 
\begin{array}{c}
p_2\circ 
\widehat{T}^{n_P}(P)>s \\ p_2\circ \widehat{T}^{n_Q}(Q)<l 
\end{array}
\right. $.
\end{theorem}

For a proof see \cite{Katok}.

As we have already said, in the proof of theorems (\ref{ImpQper}) and (\ref
{entropia}) we use some results from the Nielsen-Thurston theory of
classification of homeomorphisms of surfaces up to isotopy and some results
due to M. Handel.

The following is a brief summary of these results, taken from \cite{LMack}.
For more information and proofs see \cite{T}, \cite{FLP} and \cite{HT}.

Let $M$ be a compact, connected oriented surface possibly with boundary, and 
$f:M\rightarrow M$ be a homeomorphism. Two homeomorphisms are said to be
isotopic if they are homotopic via homeomorphisms. In fact, for closed
orientable surfaces, all homotopic pairs of homeomorphisms are isotopic \cite
{E1}.

There are two basic types of homeomorphisms which appear in the
Nielsen-Thurston classification: the finite order homeomorphisms and the
pseudo-Anosov ones.

A homeomorphism $f$ is said to be of finite order if $f^n=id$ for some $n\in 
{\rm I}\negthinspace {\rm N.}$ The least such $n$ is called the order of $f.$
Finite order homeomorphisms have topological entropy zero.

A homeomorphism $f$ is said to be pseudo-Anosov if there is a real number $%
\lambda >1$ and a pair of transverse measured foliations $F^S$ and $F^U$
such that $f(F^S)=\lambda ^{-1}F^S$ and $f(F^U)=\lambda F^U.$ Pseudo-Anosov
homeomorphisms are topologically transitive, have positive topological
entropy, and have Markov partitions \cite{FLP}.

A homeomorphism $f$ is said to be reducible by a system 
$$
C=\stackrel{n}{\stackunder{i=1}{\cup }}C_i 
$$

of disjoint simple closed curves $C_1,...,C_n$ (called reducing curves) if

(1) $\forall i,$ $C_i$ is not homotopic to a point, nor to a component of $%
\partial M,$

(2) $\forall i\neq j,$ $C_i$ is not homotopic to $C_j,$

(3) $C$ is invariant under $f.$

\begin{theorem}
\label{Thu}: If the Euler characteristic $\chi (M)<0,$ then every
homeomorphism $f:M\rightarrow M$ is isotopic to a homeomorphism $%
F:M\rightarrow M$ such that either

(a) $F$ is of finite order,

(b) $F$ is pseudo-Anosov, or

(c) $F$ is reducible by a system of curves $C.$
\end{theorem}

Homeomorphisms $F$ as in theorem (\ref{Thu}) are called Thurston canonical
forms for $f.$

\begin{theorem}
\label{EntThu}: If $f$ is pseudo-Anosov and $g$ is isotopic to $f,$ then $%
h(g)\geq h(f).$
\end{theorem}

Some results due to M. Handel can be trivially adapted to the situation
studied here. To be more precise, we can change in propositions 1.1 and 1.2
of \cite{Handel}, annulus homeomorphisms by torus homeomorphisms homotopic
to the map $LM:{\rm T^{2}\rightarrow T^{2},}$ which is the torus map induced
by the following linear map of the plane:

\begin{equation}
\label{LM}\left( 
\begin{array}{c}
\widetilde{\phi} ^{\prime } \\ \widetilde{I}^{\prime } 
\end{array}
\right) =\left( 
\begin{array}{cc}
1 & 1 \\ 
0 & 1 
\end{array}
\right) .\left( 
\begin{array}{c}
\widetilde{\phi} \\ \widetilde{I} 
\end{array}
\right) 
\end{equation}

In our case we also have to present appropriate definitions for rotation
number and rotation set.

Given a homeomorphism $\overline{f}:{\rm T^2\rightarrow T^2}$ that is
homotopic to $LM$ and a lift of $\overline{f}$ to the cylinder, $\widehat{f}%
:S^1\times {\rm I}\negthinspace {\rm R\rightarrow }S^1\times {\rm I}%
\negthinspace {\rm R,}$ we define the vertical rotation set as 
$$
\rho _V(\widehat{f})=\cup \rho _V(\widehat{f},z) 
$$
where the union is taken over all $z\in {\rm T^2}$ such that the vertical
rotation number ($\widehat{z}\in S^1\times {\rm I}\negthinspace {\rm R}$ is
any lift of $z\in {\rm T^2}$) 
$$
\rho _V(\widehat{f},z)=\stackunder{n\rightarrow \infty }{\lim }\frac{%
p_2\circ \widehat{f}^n(\widehat{z})-p_2(\widehat{z})}n 
$$
exists.

We say that $\overline{f}:{\rm T^2\rightarrow T^2}$ is pseudo-Anosov
relative to a finite invariant set $Q\subset {\rm T^2}$ if it satisfies all
of the properties of a pseudo-Anosov homeomorphism except that the
associated stable and unstable foliations may have 1-prolonged singularities
at points in $Q.$ As a last definition, for every set $A\subset {\rm T^2}$
let $\widehat{A}\subset S^1\times {\rm I}\negthinspace {\rm R}$ be the full
(cylinder) pre-image of $A{\rm .\ }$Now we are ready to state the modified
versions of propositions 1.1 and 1.2 of \cite{Handel}:

\begin{proposition}
\label{mhe1} (modified 1.1): If $\overline{f}:{\rm T^2\rightarrow T^2}$
homotopic to $LM,$ is pseudo-Anosov relative to some finite invariant set $%
Q, $ then $\rho _V(\widehat{f})$ is a closed interval. For each $\omega \in
\rho _V(\widehat{f}),$ there is a compact invariant set $E_\omega \subset 
{\rm T^2}$ such that $\rho _V(\widehat{f},z)=\omega $ for all $z\in E_\omega 
$. Moreover, if $\omega \in int\left( \rho _V(\widehat{f})\right) ,$ then we
may choose $E_\omega \subset {\rm T^2\backslash Q.}$
\end{proposition}

\begin{proof} 
As in \cite{Handel}. \hfill
\end{proof}

\begin{proposition}
\label{mhe2} (modified 1.2): Suppose that $\overline{f}:{\rm T^2\rightarrow
T^2}$ is pseudo-Anosov relative to a finite invariant set $Q$\ and that $
\overline{T}:{\rm T^2\rightarrow T^2}$ (induced by some element of $TQ$) is
homotopic to $\overline{f}$ relative to $Q.$ If $\widehat{f}:S^1\times {\rm I%
\negthinspace R\rightarrow }S^1\times {\rm I\negthinspace R}$ and $\widehat{T%
}:S^1\times {\rm I\negthinspace R\rightarrow }S^1\times {\rm I\negthinspace
R}$ are lifts that are equivariantly homotopic relative to $\widehat{Q}$,
then $\rho _V(\widehat{T})\supset int\left( \rho _V(\widehat{f})\right) .$
Moreover, for each $\omega \in int\left( \rho _V(\widehat{f})\right) ,$
there is a compact $\overline{T}-$invariant set $E_\omega \subset {\rm T^2}$
such that $\rho _V(\widehat{T},z)=\omega $ for all $z\in E_\omega .$
\end{proposition}

\begin{proof}
 Also as in \cite{Handel}. \hfill
\end{proof}

\section{Proofs}

From up to now, for simplicity, we will omit the $\symbol{126}$ in the
coordinates $(\widetilde{\phi },\widetilde{I})$ of the plane. First of all
we prove lemma (\ref{characrot}). This lemma is a trivial consequence of the
following result:

\begin{lemma}
\label{rotacoes}: Let $T\in TQ$ and $z,w_{+},w_{-}\in {\rm I}{\rm 
\negthinspace R^2}$ be points such that: 
$$
\begin{array}{l}
i) 
\text{ }\left| p_2\circ T^n(z)\right| <C,\text{ for all }n\geq 0\text{ and
some constant }C>0 \\ ii)\text{ }p_2\circ T^n(w_{\pm })\stackrel{%
n\rightarrow \infty }{\rightarrow }\pm \infty 
\end{array}
$$
So we have: 
$$
\begin{array}{l}
i) 
\text{ }\exists K>0,\text{ such that for all }n>0,\text{ }\left| \dfrac{%
p_1\circ T^n(z)-p_1(z)}n\right| <K \\ ii)\text{ }\dfrac{p_1\circ T^n(w_{\pm
})-p_1(w_{\pm })}n\stackrel{n\rightarrow \infty }{\rightarrow }\pm \infty 
\end{array}
$$
\end{lemma}

\begin{proof}

As the proofs for $w_{+}$ and $w_{-}$ are equal, we only analyze $w_{+}$,
which will be called just $w$. For all $n>0$ we define $z=(\phi _z^0,I_z^0),$
$\phi _z^n=p_1\circ T^n(z),$ $I_z^n=p_2\circ T^n(z)$ and $w=(\phi
_w^0,I_w^0),$ $\phi _w^n=p_1\circ T^n(w),$ $I_w^n=p_2\circ T^n(w)$. From the
initial hypothesis, $\left| I_z^j\right| <C$ for all $j>0,$ so defining the
following $\phi -$periodic function $\widetilde{T_\phi }(\phi ,I)=T_\phi
(\phi ,I)-\phi ,$ there is a constant $K>0$ such that:

$$
\left| \frac{p_1\circ T^n(z)-p_1(z)}n\right| \leq \frac{\stackrel{n-1}{\stackunder{j=0}{\sum }}\left| \widetilde{T_\phi }(\phi _z^j,I_z^j)\right| }n<\frac{n.K}n=K. 
$$

Now we write $I_w^n=I_{w0}^n+k_n,$ with $I_{w0}^n\in [0,1)$ and $k_n\in {\rm %
Z\negthinspace \negthinspace Z}.$ Of course the hypothesis in the lemma
implies that $k_n\stackrel{n\rightarrow \infty }{\rightarrow }\infty ,$
because $p_2\circ T^n(w)\stackrel{n\rightarrow \infty }{\rightarrow }\infty .
$ As above $\exists $ $\overline{K}>0$ such that for all $j>0,$ $\left| 
\widetilde{T_\phi }(\phi _w^j,I_{w0}^j)\right| \leq \stackunder{(\phi ,I)\in
[0,1]^2}{\max }\left| \widetilde{T_\phi }(\phi ,I)\right| <\overline{K}.$

So for all $n>0$:

$$
\frac{p_1\circ T^n(w)-p_1(w)}n=\frac{\stackrel{n-1}{\stackunder{j=0}{\sum }}
\widetilde{T_\phi }(\phi _w^j,I_{w0}^j)+\stackrel{n-1}{\stackunder{j=0}{\sum 
}}k_j}n>-\overline{K}+\frac{\stackrel{n-1}{\stackunder{j=0}{\sum }}k_j}n 
$$

In order to finish the proof that $\stackunder{n\rightarrow \infty }{\lim }
\dfrac{p_1\circ T^n(w)-p_1(w)}n=\infty ,$ we remember the Cesaro theorem,
which says that $\stackunder{j\rightarrow \infty }{\lim }k_j=\infty $ $%
\Rightarrow $ $\stackunder{n\rightarrow \infty }{\lim }\dfrac{\stackrel{n-1}{\stackunder{j=0}{\sum }}k_j}n=\infty $\hfill 
\end{proof}

\vskip0.2truecm

The following is a very important lemma.

\begin{lemma}
\label{prop}: Given $T\in TQ$ such that $\widehat{T}$ is $\mu $-exact and $T$
does not have R.I.C's we have:

For all $k\in {\rm I}\negthinspace {\rm N}^{*}$, $\exists $ $N>0$ and $
\widetilde{P}=(\phi _{\widetilde{P}},I_{\widetilde{P}})\in [0,1]^2,$ such
that $T^N(\widetilde{P})=$

$T^N(\phi _{\widetilde{P}},I_{\widetilde{P}})=(\phi _{\widetilde{P}}+s,I_{ 
\widetilde{P}}^N)$ for some $s\in {\rm Z\negthinspace \negthinspace Z},$
with $I_{\widetilde{P}}^N>I_{\widetilde{P}}+k.$
\end{lemma}

\begin{proof}

The proof will be done by contradiction. Suppose that exists $k_{0}\geq 1,$
such that $\forall $ $N>0,$ there is no $\widetilde{P}\in [0,1]^{2}$ such
that $T^{N}(\widetilde{P})=T^{N}(\phi _{\widetilde{P}},I_{\widetilde{P}})=(\phi _{\widetilde{P}}+s,I_{\widetilde{P}}^{N})$ for some $s\in {\rm %
Z\!\!Z},$ with $I_{\widetilde{P}}^{N}>1+k_{0}\geq I_{\widetilde{P}}+k_{0}.$

First let us note that given a map $T\in TQ,$ $\exists $ $a>0,$ such that $%
\forall$ $Q\in {\rm I}\negthinspace {\rm R}^{{\rm 2}}{\rm ,}$ $p_2\circ
\left(T(Q)-Q\right) >-a.$ In fact, from the definition of the set $TQ,$ we
just have to take $a>-\stackunder{Q\in [0,1]^2}{\inf }p_2\circ \left(
T(Q)-Q\right),$ because as $[0,1]^2$ is compact, $a<\infty .$

All $T\in TQ$ can be written in the following way,

$$
T:\left\{ 
\begin{array}{c}
\phi ^{\prime }=T_\phi (\phi ,I) \\ 
I^{\prime }=T_I(\phi ,I) 
\end{array}
\right. 
$$
and for all $(\phi ,I)\in {\rm I}\negthinspace {\rm R^2}$ we have the
following estimates:

\begin{equation}
\label{gf} \exists \text{ }b>0,\text{ such that }\left| \frac{\partial
T_\phi }{\partial \phi }\right| <b 
\end{equation}

\begin{equation}
\label{g} \exists \text{ }K>0,\text{ such that }\frac{\partial T_\phi }{\partial I}\geq K\text{ (twist condition)} 
\end{equation}

As $T$ does not have R.I.C's, theorem (\ref{orbclimb}) implies that

$$
\exists \text{ }P=(\phi _{P},I_{P})\in [0,1]^{2}\text{ and }N_{1}>1\text{
such that:} 
$$

$$
p_2\circ T^{N_1}(P)>\left[ (k_0+2)+\frac{(2+b)}K\right] +a 
$$

A very natural thing is to look for a point $\widetilde{P}$ as described
above, in the line segment $r=\left\{ (\phi ,I)\in [0,1]^2:\;\phi =\phi
_P\right\} .$

First, let us define 
\begin{equation}
\label{chma}Max.H.L(T^{n}(r))=\stackunder{x,y\in [0,1]}{\sup }\left|
p_{1}\circ T^{n}(\phi _{P},x)-p_{1}\circ T^{n}(\phi _{P},y)\right| . 
\end{equation}

It is clear that 
\begin{equation}
\label{chmaxde} Max.H.L(T^n(r))\geq \left| p_1\circ T^n(\phi _P,0)-p_1\circ
T^n(\phi _P,1)\right| =n\stackrel{n\rightarrow \infty }{\rightarrow }\infty
. 
\end{equation}

So, for all $n>1,$ $\exists $ at least one $s\in {\rm Z\negthinspace
\negthinspace Z},$ such that $\phi _P+s\in p_1\left( T^n(r)\right) .$

The hypothesis we want to contradict implies that for all $n>0$ and $Q\in r,$
such that

\begin{equation}
\label{hipcont1} p_1\circ T^n(Q)=\phi _P\text{ }(mod\text{ }1), 
\end{equation}
we have:

\begin{equation}
\label{hipcont2} p_2\circ T^n(Q)\leq (k_0+1) 
\end{equation}

As $p_2\circ T^{N_1}(P)>\left[ (k_0+2)+\frac{(2+b)}K\right] +a,$ $\exists $ $%
P_1\in r$ such that:

$$
\begin{array}{c}
p_2\circ T^{N_1}(P_1)=(k_0+2)+a \\ 
\text{and} \\ 
\begin{array}{l}
\forall 
\text{ }Q\in \text{ }\overline{PP_1}\subset r, \\ p_2\circ T^{N_1}(Q)\geq
(k_0+2)+a 
\end{array}
\end{array}
$$
The reason why such a point $P_1$ exists is the following: As $N_1>1,$ $%
\exists $ at least one $s\in {\rm Z\negthinspace \negthinspace Z}$ such that 
$\phi _P+s\in p_1\left( T^{N_1}(r)\right) .$ Thus, from (\ref{hipcont1}) and
(\ref{hipcont2}), $T^{N_1}(r)$ must cross the line $l$ given by: $l=\left\{
(\phi ,(k_0+2)+a),\text{ with }\phi \in {\rm I}\negthinspace {\rm R}\right\} 
$

Also from (\ref{hipcont1}) and (\ref{hipcont2}) we have that: 
$$
\stackunder{Q,R\in \overline{PP_1}}{\sup }\left| p_1\circ
T^{N_1}(R)-p_1\circ T^{N_1}(Q)\right| <1 
$$

Now let $\gamma _{N_1}:J\rightarrow {\rm I}\negthinspace {\rm R^2}$ be the
following curve: 
\begin{equation}
\label{Jdef}\gamma _{N_1}(t)=T^{N_1}(\phi _P,t),\text{ }t\in J=\text{interval whose extremes are }I_P\text{ and }I_{P_1} 
\end{equation}
It is clear that it satisfies the following inequalities: 
$$
\begin{array}{c}
p_2\circ \gamma _{N_1}(I_P)-p_2\circ \gamma _{N_1}(I_{P_1})> 
\frac{(2+b)}K \\ \text{ } \\ \text{ }\stackunder{t,s\in J}{\sup }\left|
p_1\circ \gamma _{N_1}(t)-p_1\circ \gamma _{N_1}(s)\right| <1 
\end{array}
$$

\begin{description}
\item[Claim 1]  : Given a continuous curve $\gamma :J=[\alpha ,\beta
]\rightarrow {\rm I}\negthinspace 
{\rm R^2,}$ with 
\begin{equation}
\label{gira1}\stackunder{t,s\in J}{\sup }\left| p_1\circ \gamma (t)-p_1\circ
\gamma (s)\right| <1
\end{equation}

\begin{equation}
\label{gira2}\left| p_2\circ \gamma (\beta )-p_2\circ \gamma (\alpha
)\right| >\frac{(2+b)}K
\end{equation}

Then $\exists $ $s\in {\rm Z\negthinspace \negthinspace Z},$ such that $\phi
_P+s\in p_1\left( T\circ \gamma (J)\right) .$
\end{description}

\begin{proof}

$$
\begin{array}{l}
\stackunder{t,s\in J}{\sup }\left| p_1\circ T\circ \gamma (t)-p_1\circ
T\circ \gamma (s)\right| = \\ =
\stackunder{t,s\in J}{\sup }\left| T_\phi \circ \gamma (t)-T_\phi \circ
\gamma (s)\right| \geq \left| T_\phi \circ \gamma (\beta )-T_\phi \circ
\gamma (\alpha )\right| \geq  \\ \geq -b+K.\frac{(2+b)}K=2
\end{array}
$$
So the claim is proved.\hfill
\end{proof}

\vskip0.2truecm

$\gamma _{N_1}(t)$ (see (\ref{Jdef})) satisfies the claim hypothesis, by
construction. So $\exists $ $s\in {\rm Z\negthinspace 
\negthinspace Z}$ such that $\phi _P+s\in p_1\left( T\circ \gamma
_{N_1}(J)\right) =p_1\left( T^{N_1+1}(\overline{PP_1})\right) $.

As $\stackunder{t\in J}{\inf }$ $p_2\left( \gamma _{N_1}(t)\right)
=p_2\left( \gamma _{N_1}(I_{P_1})\right) =(k_0+2)+a,$ from the choice of $%
a>0 $ we get that $\stackunder{t\in J}{\inf }$ $p_2\left( T\circ \gamma
_{N_1}(t)\right) >(k_0+2).$

So there is $\overline{t}\in J$ and $\overline{P}=(\phi _P,\overline{t})\in
r $ such that:

$$
\begin{array}{l}
p_1\circ T^{N_1+1}( 
\overline{P})=\phi _P\text{ }(mod\text{ }1) \\ p_2\circ T^{N_1+1}(\overline{P})>(k_0+2) 
\end{array}
$$

And this contradicts (\ref{hipcont1}) and (\ref{hipcont2}). So for all $%
k\geq 1$, $\exists $ $N>0$ and $\widetilde{P}\in r,$ such that $T^{N}( 
\widetilde{P})=T^{N}(\phi _{P},I_{\widetilde{P}})=(\phi _{P}+s,I_{\widetilde{P}}^{N})$ for some $s\in {\rm Z\!\!Z},$ with $I_{\widetilde{P}}^{N}>I_{ 
\widetilde{P}}+k $. \hfill 
\end{proof}

\vskip0.2truecm

Remark:

$\bullet $ Of course for all $k\leq -1,$ $k\in {\rm Z\!\!Z}$, there are also 
$N>0$ and $\widetilde{Q}=(\phi _{\widetilde{Q}},I_{\widetilde{Q}})\in
[0,1]^{2}, $ such that $T^{N}(\widetilde{Q})=T^{N}(\phi _{\widetilde{Q}},I_{ 
\widetilde{Q}})=(\phi _{\widetilde{Q}}+s,I_{\widetilde{Q}}^{N})$ for some $%
s\in {\rm Z\!\!Z},$ with $I_{\widetilde{Q}}^{N}<I_{\widetilde{Q}}+(k-1).$
The proof in this case is completely similar to the above one, because as $T$
does not have R.I.C's, again by theorem (\ref{orbclimb}) for all $l<0$ there
exist $Q=(\phi _{Q},I_{Q})\in [0,1]^{2}$ and $n_{Q}>1$ such that $p_{2}\circ
T^{n_{Q}}(Q)<l.$

Below we prove the main results of this paper.

\vskip0.2truecm

\begin{proof}{\it of theorem (\ref{PtosPerQuoc})}

As the 2 cases, $k>0$ and $k<0$ are completely similar, let us fix $k>0.$

$\left( \Rightarrow \right) $

If $\overline{T}$ has a periodic point $P,$ with $\rho _V(P)=\frac kN,$ for
some $k>0$ and $N>0,$ then $p_2\circ \widehat{T}^n(P)\stackrel{n\rightarrow
\pm \infty }{\rightarrow }\pm \infty ,$ which implies that there can be no
R.I.C.

$\left( \Leftarrow \right) $

To prove the existence of a periodic orbit with $\rho _V=\frac kN,$ for a
given $k>0$ and some $N>0$ sufficiently large, it is enough to show that
there exists a point $P\in S^1\times {\rm I}\negthinspace {\rm R}$ such that:

\begin{equation}
\label{Pij1}\widehat{T}^N(P)=P+(0,k) 
\end{equation}

As $T\in TQ,$ for each $(s,l)\in {\rm Z\negthinspace \negthinspace Z}^2$ and 
$N>0$ the sets $C(s,N)$, defined in lemma (\ref{Lcal1}), satisfy: $%
C(s+l.N,N)=C(s,N)+(0,l)$

So, for each fixed $N>0$, there are only $N$ distinct sets of this type: $%
C(0,N),C(1,N),...,C(N-1,N)$

The others are just integer vertical translations of them. Another trivial
remark about the sets $C(s,N)$ is: $C(s,N)\cap C(r,N)=\emptyset ,$ if $s\neq
r$

For all $s\in {\rm Z\negthinspace \negthinspace Z,}$ we get from lemma (\ref
{inter}) that $\widehat{T}(C(s,N))\cap C(s,N)\neq \emptyset .$ So we can
apply lemma (\ref{allrotnum}) and conclude that $\exists $ $\overline{P}%
_s\in C(s,N)$ such that $\widehat{T}^N(\overline{P}_s)=\overline{P}_s.$

From lemma (\ref{prop}), for the given $k>0$, $\exists $ $N>0$ and $
\widetilde{P}=(\phi _{\widetilde{P}},I_{\widetilde{P}})\in S^1\times {\rm I}%
\negthinspace {\rm R},$ such that $\widehat{T}^N(\widetilde{P})=\widehat{T}%
^N(\phi _{\widetilde{P}},I_{\widetilde{P}})=(\phi _{\widetilde{P}},I_{
\widetilde{P}}^N),$ with $I_{\widetilde{P}}^N>I_{\widetilde{P}}+k.$ So $
\widetilde{P}\in K(\widetilde{s},N)$ (see expression (\ref{Kpq}))$,$ for a
certain $\widetilde{s}\in {\rm Z\negthinspace \negthinspace Z}$ and $%
p_2\circ \widehat{T}^N(\widetilde{P})-p_2(\widetilde{P})>k.$ As $\widetilde{P}\in K(\widetilde{s},N),$ we get that $\mu ^{-}(\phi _{\widetilde{P}})\leq
p_2(\widetilde{P})$ and $\nu $$^{+}(\phi _{\widetilde{P}})\geq p_2\circ 
\widehat{T}^N(\widetilde{P}),$ which implies that $\nu ^{+}(\phi _{
\widetilde{P}})-\mu ^{-}(\phi _{\widetilde{P}})\geq p_2\circ \widehat{T}^N(
\widetilde{P})-p_2(\widetilde{P})>k.$

From corollary (\ref{ofpre}) we get that $\widehat{T}^N(\phi _{\widetilde{P}},\mu ^{-}(\phi _{\widetilde{P}}))=(\phi _{\widetilde{P}},\nu ^{+}(\phi _{
\widetilde{P}})),$ so defining $\widetilde{\widetilde{P}}=(\phi _{\widetilde{P}},\mu ^{-}(\phi _{\widetilde{P}}))\in C(\widetilde{s},N)$ (see lemma (\ref
{graphu})) we have: $p_2\circ \widehat{T}^N(\widetilde{\widetilde{P}})-p_2(
\widetilde{\widetilde{P}})>k.$

And as we proved above, $\exists $ $\overline{P}_{\widetilde{s}}\in C(
\widetilde{s},N)$ such that $p_2\circ \widehat{T}^N(\overline{P}_{\widetilde{s}})-p_2(\overline{P}_{\widetilde{s}})=0.$

So as $C(\widetilde{s},N)$ is connected, $\exists $ $P\in C(\widetilde{s},N)$
such that:

$$
p_2\circ \widehat{T}^N(P)=p_2(P)+k 
$$

And the theorem is proved. Now we present an alternative proof, suggested by
a referee, which is much shorter. We decided to mantain the original proof
because it is based on lemma (\ref{prop}), which will be used in future
works, so we wanted to keep it in the present paper.

For a given $k>0,$ we are going to prove the existence of a point $P\in C(0,
\overline{N})$ such that $\widehat{T}^{\overline{N}}(P)=P+(0,k),$ for a
sufficiently large $\overline{N}.$ As $\widehat{T}$ is $\mu $-exact, we get
that there is a point $P_0=(\phi _0,I_0)\in C(0,1)$ such that $\widehat{T}%
(P_0)=P_0.$ For any given $N>0,$ let $\mu ^{-},\mu ^{+},\nu ^{-},\nu ^{+}$
be the maps associated to $C(0,N).$ From the choice of $P_0$ we get that $%
\mu ^{-}(\phi _0)\leq I_0\leq \nu ^{+}(\phi _0).$ In the proof of lemma (\ref
{Lcal2}) (see \cite{lecalvez1}), the following property for the lift of the
map $\mu ^{-}$ to ${\rm I\negthinspace R}$ is obtained (we are denoting the
lift also by $\mu ^{-}$): for any $n\in \{1,2,...,N\}$, the point $T^n(\phi
,\mu ^{-}(\phi ))$ is the first point where the image of $\phi \times {\rm I%
\negthinspace R}$ by $T^n$ meets the vertical passing through $T^n(\phi ,\mu
^{-}(\phi )),$ and for the same reasons, we get that $p_1\circ T^n(\phi
^{\prime },\mu ^{-}(\phi ^{\prime }))<p_1\circ T^n(\phi ,\mu ^{-}(\phi ))$
if $\phi ^{\prime }<\phi .$ So the image by $T^n$ of the graph of $\mu ^{-}$
is also a graph and the order given by $p_1$ is preserved. Moreover, as $T$
is a twist map, we can prove that (see lemma 13.1.1, page 424 of \cite{Katok}%
) if $\phi >\phi ^{\prime },$ then $\mu ^{-}(\phi )-\mu ^{-}(\phi ^{\prime
})\geq -\text{cotan}(\beta ).(\phi -\phi ^{\prime }),$ where $\beta $ is a
uniform angle of deviation for $T$. By periodicity of $\mu ^{-}$ we get that 
$\max $ $\mu ^{-}-\min $ $\mu ^{-}\leq \text{cotan}(\beta )$ and analogous
inequalities for the other maps.

As in the above proof, we know that $p_2\circ \widehat{T}^N-p_2$ vanishes on 
$C(0,N)$. Suppose that this map does not take the value $k$ on $C(0,N).$
Then as $C(0,N)$ is compact, it is sctrictly smaller and we have $\nu
^{+}(\phi )-\mu ^{-}(\phi )<k,$ for all $\phi \in S^1.$\ So for any $\phi
\in S^1,$ we get the following estimates: 
$$
\begin{array}{l}
\mu ^{-}(\phi )=\mu ^{-}(\phi )-\mu ^{-}(\phi _0)+\mu ^{-}(\phi _0)-\nu
^{+}(\phi _0)+\nu ^{+}(\phi _0)>-
\text{cotan}(\beta )-k+I_0 \\ \nu ^{+}(\phi )=\nu ^{+}(\phi )-\nu ^{+}(\phi
_0)+\nu ^{+}(\phi _0)-\mu ^{-}(\phi _0)+\mu ^{-}(\phi _0)<\text{cotan}(\beta
)+k+I_0 
\end{array}
$$

And the above inequalities imply that 
$$
\widehat{T}^N(S^1\times ]-\infty ,-\text{cotan}(\beta )-k+I_0])\subset
S^1\times ]-\infty ,\text{cotan}(\beta )+k+I_0], 
$$
which can not hold for all $N>0$ by theorem (\ref{orbclimb}). \hfill 
\end{proof}

\vskip0.2truecm

\begin{proof}{\it of theorem (\ref{ImplicPeriod})}

Again we fix $k>0\Rightarrow k^{\prime }>0.$ The case $k<0$ is completely
similar.

By contradiction, suppose that for some $0<\dfrac{k^{\prime }}{N^{\prime }}%
<\dfrac kN$ and any fixed $s\in {\rm Z\negthinspace \negthinspace Z:}$ 
$$
p_2\circ \widehat{T}^{N^{\prime }}(Q)-p_2(Q)-k^{\prime }\leq 0,\text{ }%
\forall Q\in C(s,N^{\prime }). 
$$

So in particular, we have: $\nu ^{+}(\phi )-\mu ^{-}(\phi )-k^{\prime }\leq
0 $ for all $\phi \in S^1.$ This means that the unbounded connected
component of $C(s,N^{\prime })^c,$ which is below $C(s,N^{\prime })$ and we
denote by $U,$ satisfies the following equation: $\widehat{T}^{N^{\prime
}}(U)-(0,k^{\prime })\subset U,$ so $\widehat{T}^{i.N^{\prime
}}(U)-(0,i.k^{\prime })\subset U,$ for all $i>0.$ Now let us choose a point $%
P\in U,$ such that 
\begin{equation}
\label{ksobn}\stackunder{n\rightarrow \infty }{\lim }\frac{p_2\circ \widehat{%
T}^n(P)-p_2(P)}n=\frac kN 
\end{equation}

So we get that for all $i>0,$ $\left[ p_2\circ \widehat{T}^{i.N^{\prime
}}(P)-p_2(P)-i.k^{\prime }\right] \leq C-p_2(P),$ where $C=\sup
\{p_2(x):x\in C(s,N^{\prime })\}.$ And this implies that:%
$$
\stackunder{i\rightarrow \infty }{\lim }\frac{p_2\circ \widehat{T}%
^{i.N^{\prime }}(P)-p_2(P)}{i.N^{\prime }}\leq \frac{k^{\prime }}{N^{\prime }%
}\Rightarrow \stackunder{n\rightarrow \infty }{\lim }\frac{p_2\circ \widehat{%
T}^n(P)-p_2(P)}n\leq \frac{k^{\prime }}{N^{\prime }}, 
$$
which contradicts (\ref{ksobn}). So there is a point $P_1=(\phi _1,I_1)\in
C(s,N^{\prime })$ such that $p_2\circ \widehat{T}^{N^{\prime
}}(P_1)-p_2(P_1)-k^{\prime }>0$ and from the $\mu $-exactness of $\widehat{T}%
,$ $\exists $ $P_0=(\phi _0,I_0)\in C(s,N^{\prime })$ such that $p_2\circ 
\widehat{T}^{N^{\prime }}(P_0)-p_2(P_0)<0.$ Now let $\Delta _0$ and $\Delta
_1$ be the proper simple arcs given by:

$$
\begin{array}{c}
\Delta _0=\{\phi _0\}\times [\mu ^{+}(\phi _0),+\infty [
\text{ }\cup \text{ }\widehat{T}^{-N^{\prime }}(\{\phi _0\}\times ]-\infty
,\nu ^{-}(\phi _0)] \\ \Delta _1=\{\phi _1\}\times ]-\infty ,\mu ^{-}(\phi
_1)]\text{ }\cup \text{ }\widehat{T}^{-N^{\prime }}(\{\phi _1\}\times [\nu
^{+}(\phi _1),+\infty [
\end{array}
$$
It is easy to see that $\Delta _0\cap C(s,N^{\prime })=(\phi _0,\mu
^{+}(\phi _0)),$ $\Delta _1\cap C(s,N^{\prime })=(\phi _1,\mu ^{-}(\phi _1))$
and that $(\Delta _0\cup \Delta _1)^c$ is an open set that divides $%
C(s,N^{\prime })$ into 2 connected components, $C_1$ and $C_2$ ($%
C(s,N^{\prime })=C_1\cup C_2$), such that $C_1\cap C_2=(\phi _0,\mu
^{+}(\phi _0))\cup (\phi _1,\mu ^{-}(\phi _1)).$ Therefore the function $%
p_2\circ \widehat{T}^{N^{\prime }}-p_2-k^{\prime }$ has at least one zero in
each $C_i.$ \hfill 
\end{proof}

\vskip0.2truecm

\begin{proof}{\it of theorem (\ref{ImpQper})}

The proof will be divided into 2 cases\ (as before we fix $\omega
>0\Rightarrow \omega ^{\prime }>0$):

Case I)\ $\omega \in $ {\footnotesize {\sf l}}${\rm \negthinspace
\negthinspace }{\sf Q}$.

As $\omega ^{\prime }\in {\rm I}\negthinspace {\rm R\backslash }$%
{\footnotesize {\sf l}}${\rm \negthinspace \negthinspace }{\sf Q}$, there is
a sequence 
$$
\frac{p_i}{q_i}\stackrel{i\rightarrow \infty }{\rightarrow }\omega ^{\prime
},\text{ with }0<\frac{p_i}{q_i}<\omega ,\text{ }\forall i>0 
$$

and (from theorem (\ref{ImplicPeriod})) a family of periodic orbits

$$
E_i=\{P_1^i,P_2^i,...,P_{q_i}^i\}\subset {\rm T^2,}\text{ with }\rho
_V(E_i)= \frac{p_i}{q_i}\text{ .} 
$$

So in the Hausdorff topology there is a subsequence $E_{i_n}\stackrel{n\rightarrow \infty }{\rightarrow }E\subset {\rm T^2}$ that for simplicity
we will call $E_n.$ The convergence in the Hausdorff topology means that:
Given $\epsilon >0,$ $\exists $ $n_0\in {\rm I\negthinspace N,}$ such that
for all $n\geq n_0,$ $E_n\subset B_\epsilon (E)$ and $E\subset B_\epsilon
(E_n),$ where $B_\epsilon (\bullet )$ is the $\epsilon $ neighborhood of the
given set$.$

In this way, for all $z\in E$ there is a sequence $z_n\stackrel{n\rightarrow
\infty }{\rightarrow }z,$ such that $z_n\in E_n.$ But there is still a
problem to obtain that for all $z\in E,$ $\rho _V(z)=\omega ^{\prime }.$
Because we do not have any control under the uniformity of the vertical
rotation numbers of the family of orbits $E_n.$ In the Aubry-Mather case,
the periodic orbits, whose limit in the Hausdorff topology is a
quasi-periodic set, have a very strong uniformity condition; they are of
Birkhoff type (see for instance \cite{katok}). Indeed, if we knew that given 
$\epsilon >0,$ $\exists $ $\widehat{i}(\epsilon )>0$ ($\widehat{i}$
independent of $n$)$,$ such that $\forall $ $n>0$ and $\forall $ $z_n\in E_n$
$$
\left| \frac{p_2\circ T^i(z_n)-p_2(z_n)}i-\frac{p_n}{q_n}\right| <\epsilon , 
\text{ for all}\ i>\widehat{i}, 
$$
then the problem would be solved. In order to overcome this problem we use
the following important lemma that is a consequence of propositions (\ref
{mhe1}), (\ref{mhe2}), some ideas from \cite{LMack} and some results from
the Nielsen-Thurston theory:

\begin{lemma}
\label{adapt}: Under the hypothesis of theorem (\ref{ImpQper}), for all $%
{\rm \omega ^{\prime }\in (0,\omega )\backslash }${\footnotesize {\sf l}}$%
{\rm \negthinspace \negthinspace }{\sf Q}$, there is a quasi-periodic set $
\overline{E},$ such that $\rho _V(\overline{E})=\omega ^{\prime }.$
\end{lemma}

\begin{proof}
See the end of section 3. \hfill
\end{proof}

Case II)\ $\omega \notin $ {\footnotesize {\sf l}}${\rm \negthinspace
\negthinspace }{\sf Q}$.

From case I), we just have to prove that for all $0<\frac pq<\omega $ there
is a $q$-periodic orbit with $\rho _V=\frac pq.$ The proof of this fact is
identical to the proof of theorem (\ref{ImplicPeriod}), so we omit it. 
\hfill 
\end{proof}

\vskip0.2truecm

We still have to prove lemma (\ref{adapt}) and theorem (\ref{entropia}). The
following are auxiliary results that are important in these proofs.

\begin{lemma}
\label{curvinv}: Let $\overline{T}:{\rm T^2\rightarrow T^2}$ (${\rm T^2=I%
\negthinspace R^2/Z\negthinspace \negthinspace Z^2}$) be a homeomorphism
homotopic to $LM$ (see (\ref{LM})), and let $C\subset {\rm T^2}$ be a
homotopically non-trivial simple closed curve, $\overline{T}^s-$invariant,
for some $s>0.$ Then $C$ is a rotational simple closed curve on the cylinder 
$S^1\times {\rm I\negthinspace R=}${\rm (}${\rm I}\negthinspace 
{\rm R/Z\negthinspace \negthinspace Z}${\rm )}$\times {\rm I}\negthinspace 
{\rm R}$. Moreover, $[C]$ (homotopy class of $C$) is the only homotopy class
of simple closed curves on the torus that is preserved by iterates of $
\overline{T}$.
\end{lemma}

\begin{proof}

The action of $\overline{T}$ on $\pi _1{\rm (T^2)}$ is given by 
$$
\overline{T}_{*}([C])=\overline{T}_{*}(c_\phi ,c_I)=\left( 
\begin{array}{cc}
1 & 1 \\ 
0 & 1
\end{array}
\right) \left( 
\begin{array}{c}
c_\phi  \\ 
c_I
\end{array}
\right) , 
$$
and the eigenvector corresponding to the eigenvalue $1$ is $(1,0).$\hfill 
\end{proof}

\vskip0.2truecm

\begin{lemma}
\label{circcilin}: Let $\overline{f}:{\rm T^2\rightarrow T^2}$ be a
homeomorphism isotopic to $LM.$ If $\exists $ $l>0$ such that $\overline{f}^l
$ has a rotational invariant curve $\gamma $ with $[\gamma ]=x_\phi =(1,0),$
then $\overline{f}$ can not have a periodic orbit with vertical rotation
number $\rho _V\neq 0$ and another with $\rho _V=0$.
\end{lemma}

\begin{proof}

Let $\widehat{F}_0$ be a lift of $\overline{f}^l$ to the cylinder which
fixes $\widehat{\gamma },$ a lift of $\gamma \subset {\rm T^2.}$ This
implies that the vertical rotation number for $\widehat{F}_0$ of every point
is zero. So, given any lift $\widehat{F}$ of $\overline{f}^l,$ the vertical
rotation number of every point is equal. In particular, given a lift $
\widehat{f}$ of $\overline{f},$ the vertical rotation number of every point
is the same, which is what we wanted to prove. \hfill 
\end{proof}

\vskip0.2truecm

We have already seen that all $T\in TQ$ such that $\widehat{T}$ is $\mu $%
-exact and $T$ does not have R.I.C's induces a map $\overline{T}$ defined on
the torus that has periodic orbits with non-zero vertical rotation numbers.
Suppose that $\overline{T}$ has a periodic orbit with $\rho _V=\omega >0.$
Given $\omega ^{\prime }\in {\rm I\negthinspace
R\backslash }$l${\rm \negthinspace \negthinspace }{\sf Q}$, with $0<\omega
^{\prime }<\omega ,$ let us choose irreducible fractions $\frac{a_1}{b_1}$
and $\frac{a_2}{b_2},$ such that

$$
0<\frac{a_1}{b_1}<\omega ^{\prime }<\frac{a_2}{b_2}\leq \omega , 
$$

and periodic orbits $Q_1$ and $Q_2$ with $\rho _V(Q_i)=\frac{a_i}{b_i}$ and $%
\#\{Q_i\}=b_i$, for $i=1,2$ (this is possible by theorem (\ref{ImplicPeriod}%
)).

As $T\in TQ$ and $\widehat{T}$ is $\mu $-exact, from lemmas (\ref{inter})
and (\ref{allrotnum}) it is clear that $\exists $ $R\in {\rm T^2}$ such that 
$\overline{T}(R)=R$ and

$$
\begin{array}{l}
p_2\circ T( 
\widetilde{R})=p_2(\widetilde{R}) \\ p_1\circ T(\widetilde{R})=p_1( 
\widetilde{R})\text{ }(mod\text{ }1) 
\end{array}
,\text{ for any }\widetilde{R}\in p^{-1}(R). 
$$

Let $Q=Q_1\cup Q_2\cup R.$ Now we blow-up each $x\in Q$ to a circle $S_x.$
Let ${\rm T_Q^2}$ be the compact manifold (with boundary) thereby obtained ; 
${\rm T_Q^2\ }$is the compactification of ${\rm T^2\backslash Q,}$ where $%
S_x $ is a boundary component where $x$ was deleted. Now we extend $
\overline{T}:{\rm T^2\backslash Q\rightarrow T^2\backslash Q}$ to $\overline{%
T}_Q:{\rm T_Q^2\rightarrow T_Q^2}$ by defining $\overline{T}%
_Q:S_x\rightarrow S_x$ via the derivative; we just have to think of $S_x$ as
the unit circle in $T_x{\rm T^2}$ and define

$$
\overline{T}_Q(v)=\frac{D\overline{T}_x(v)}{\left\| D\overline{T}%
_x(v)\right\| },\text{ for }v\in S_x. 
$$

$\overline{T}_Q$ is continuous on ${\rm T_Q^2}$ because $\overline{T}$ is $%
C^1$ on ${\rm T^2.}$ Let $b:{\rm T_Q^2\rightarrow T^2}$ be the map that
collapses each $S_x$ onto $x.$ Then $\overline{T}\circ b=b\circ \overline{T}%
_Q.$ This gives $h(\overline{T}_Q)\geq h(\overline{T})$ (see \cite{Katok},
page 111). Actually $h(\overline{T}_Q)=h(\overline{T}),$ because each fibre $%
b^{-1}(y)$ is a simple point or an $S_x$ and the entropy of $\overline{T}$
on any of these fibres is 0 (the map on the circle induced from any linear
map has entropy 0). This construction is due to Bowen (see \cite{Bowen}).

Now we have the following:

\begin{theorem}
\label{p-anosov}: The map $\overline{T}_Q:{\rm T_Q^2\rightarrow T_Q^2}$ is
isotopic to a pseudo-Anosov homeomorphism of ${\rm T_Q^2}$.
\end{theorem}

\begin{proof}

By theorem (\ref{Thu}), $\overline{T}_Q$ is isotopic to a homeomorphism $F_Q:%
{\rm T_Q^2\rightarrow T_Q^2}$ (Thurston canonical form) such that either: 

i) $F_Q$ has finite order,

ii) $F_Q$ is pseudo-Anosov,

iii) $F_Q$ is reducible by a system of curves $C.$

We must think of ${\rm T_Q^2}$ as a torus with round disks removed, all of
the same size, each one centered in a point $x\in Q$. Let $F:{\rm %
T^2\rightarrow T^2}$ be the completion of $F_Q$, i.e. the homeomorphism
obtained by radially extending $F_Q$ into all the holes (see \cite{E2}).

It is easy to see that $F_Q$ does not have finite order, because there are
points with different rotation numbers (by construction of $Q$).

We say that a simple closed curve $\gamma $ on a torus with holes is
rotational if after filling in the holes, $\gamma $ is homotopically
non-trivial. Suppose that $F_Q$ has a rotational reducing curve $\gamma $
and let $[\gamma ]\in \pi _1({\rm T^2})$ be its homotopy class in the torus
without holes. Then, for some $n>0,$ we have: 
\[
F_Q^n(\gamma )=\gamma \Rightarrow F^n(\gamma )=\gamma . 
\]

And, as $F_{Q}$ is isotopic to $\overline{T}_{Q}$, $F$ is isotopic to $LM$.
In this way, from lemma (\ref{curvinv}) the homotopy class of $\gamma $ in
the torus ${\rm T^{2}=S^{1}\times S^{1}}$ is $[\gamma ]=x_{\phi }=(1,0)$ ($%
\gamma $ is a rotational simple closed curve in the cylinder $S^{1}\times 
{\rm I\!R}$). So, from the existence of the periodic orbits $Q_{1}$ (or $%
Q_{2}$) and $R,$ applying lemma (\ref{circcilin}) we conclude that, $F$ and
thus $F_{Q}$ do not have any rotational reducing curves.

And if $\gamma $ is a non-rotational reducing curve, then $\gamma $ must
surround at least 2 holes (because $\gamma $ is not homotopic to a component
of $\partial {\rm T_Q^2}$)${\rm .}$ These holes must have the same rotation
number and this is impossible, because $\rho _V(Q_1)\neq \rho _V(Q_2)\neq
\rho _V(R)=0$ and 2 points from the same orbit can not be surrounded by the
same curve ( by construction of $Q_1$ and $Q_2$).

So $F_Q:{\rm T_Q^2\rightarrow T_Q^2}$ is a pseudo-Anosov homeomorphism.
\hfill
\end{proof}

Now we prove theorem (\ref{entropia}):

\begin{proof}{\it of theorem (\ref{entropia}) }

We just have to see that after all the previous work, the map $\overline{T}%
_Q:{\rm T_Q^2\rightarrow T_Q^2}$ is isotopic to a pseudo-Anosov
homeomorphism of ${\rm T_Q^2}$. Then:

$h(\overline{T})=h(\overline{T}_Q)$ and $h(\overline{T}_Q)>0,$ by theorem (%
\ref{EntThu}). \hfill
\end{proof}

And finally we prove lemma (\ref{adapt}):

\begin{proof}{\it of lemma (\ref{adapt})} 

By theorem (\ref{p-anosov}), $\overline{T}_Q:{\rm T_Q^2\rightarrow T_Q^2}$
is isotopic to a pseudo-Anosov homeomorphism, $F_Q:{\rm T_Q^2\rightarrow
T_Q^2}$. So, as $Q$ is an invariant and finite set we just have to apply
propositions (\ref{mhe1}) and (\ref{mhe2}). \hfill
\end{proof}


\section{Examples and applications}

We conclude by giving some examples.

\vskip0.3truecm

1) It is obvious that the well-known Standard map $S_M:{\rm T^2\rightarrow
T^2}$ given by 
$$
S_M:\left\{ 
\begin{array}{l}
\phi ^{\prime }=\phi +I^{\prime } 
\text{ }(\func{mod}\text{ }1) \\ I^{\prime }=I-\frac k{2\pi }\sin (2\pi \phi
)\text{ }(\func{mod}\text{ }1) 
\end{array}
\text{ is induced by an element of $TQ.$}\right. 
$$
Also, it is easy to see that its generating function is: 
$$
h_{S_M}(\phi ,\phi ^{\prime })=\frac{(\phi ^{\prime }-\phi )^2}2+\frac
k{4\pi ^2}.\cos (2\pi \phi )\Rightarrow h_{S_M}(\phi +1,\phi ^{\prime
}+1)=h_{S_M}(\phi ,\phi ^{\prime }), 
$$
so $S_M$ is an exact map. In this way, as we know that for sufficiently
large $k>0,$ $S_M$ does not have R.I.C's, we can apply our previous results
to this family of maps. In fact, theorems (\ref{PtosPerQuoc}) and (\ref
{ImplicPeriod}) can be used to produce a new criteria to obtain estimates
for the parameter value $k_{cr},$ which is defined in the following way: if $%
k>k_{cr},$ then $S_M$ does not have R.I.C and for $k\leq k_{cr}$ there is at
least one R.I.C. This happens because for each $\frac 1n,$ $n\in {\rm I}%
\negthinspace {\rm N^{*},}$ there is a number $k_n,$ such that for $k\geq
k_n,$ $S_M$ has a $n$-periodic orbit with $\rho _V=\frac 1n$ and for $k<k_n$
it does not have such an orbit. From the theorems cited above, if $n>m$ then 
$k_n\leq k_m$ and $\stackunder{n\rightarrow \infty }{\lim }k_n=k_{cr}.$ In a
future work, we will try to obtain estimates for $k_{cr}$ using this method.

\vskip0.3truecm

2) In \cite{mythesis} the dynamics near a homoclinic loop to a saddle-center
equilibrium of a 2-degrees of freedom Hamiltonian system was studied by
means of an approximation of a certain Poincar\'e map.\ In an appropriate
coordinate system this map is given by: 
$$
\widehat{F}:S^1\times ]0,+\infty [\rightarrow S^1\times ]0,+\infty [,\text{
where }\widehat{F}:\left\{ 
\begin{array}{l}
\phi ^{\prime }=\mu \;(\phi )+\gamma \log (I^{\prime }) 
\text{ }(\func{mod}\text{ }\pi ) \\ I^{\prime }=J(\phi )I 
\end{array}
\right. 
$$
and 
$$
\begin{array}{l}
J(\phi )=\alpha ^2\cos {}^2(\phi )+\alpha ^{-2}\sin {}^2(\phi ) \\ 
\mu (\phi )=\arctan (\frac{\tan (\phi )}{\alpha ^2}),\text{ }\mu (0)=0 
\end{array}
. 
$$

So $J(\phi )$ is $\pi $-periodic and $\mu (\phi +\pi )=\mu (\phi )+\pi .$ In
this case $S^{1}$ will be identified with ${\rm I}\!{\rm R/(\pi Z\!\!Z).}$

A direct calculation shows that: 
$$
h_{\widehat{F}}(\phi ,\phi ^{\prime })=\gamma \exp (\frac{\phi ^{\prime
}-\mu \;(\phi )}{\gamma }) 
$$

And so 
$$
h_{\widehat{F}}(\phi +\pi ,\phi ^{\prime }+\pi )=h_{\widehat{F}}(\phi ,\phi
^{\prime }),\text{ because }\mu (\phi +\pi )=\mu (\phi )+\pi .\text{ } 
$$

Thus $\widehat{F}$ is also exact. Applying the following coordinate change 
$$
\left\{ 
\begin{array}{l}
\widetilde{\phi }=\phi \\ \widetilde{I}=\gamma \log (I) 
\end{array}
\right. 
$$
we get (omitting the $\symbol{126}$): 
$$
\widehat{F}:S^{1}\times {{\rm I}\!{\rm R\hookleftarrow }}:\left\{ 
\begin{array}{l}
\phi ^{\prime }=F_{\phi }(I,\phi )=\mu (\phi )+I^{\prime } 
\text{ }(\func{mod}\text{ }\pi ) \\ I^{\prime }=F_{I}(I,\phi )=\gamma \log
(J(\phi ))+I 
\end{array}
\right. 
$$

It is obvious that in these coordinates $\widehat{F}$ is $\mu $-exact and $%
\mu $ is given by: 
$$
\mu (A)=\stackunder{A}{\int }e^{\frac{I}{\gamma }}d\phi dI 
$$

It is also easy to see that $\widehat{F}$ induces a map $\overline{F}:{\rm %
T^2\rightarrow T^2}$ (${\rm T^2=I}\negthinspace {\rm R^2/(\pi Z%
\negthinspace
\negthinspace Z)^2}$) given by 
$$
\overline{F}:\left\{ 
\begin{array}{l}
\phi ^{\prime }= 
\overline{F}_\phi (I,\phi )=\mu (\phi )+I^{\prime }\text{ }(\func{mod}\text{ 
}\pi ) \\ I^{\prime }=\overline{F}_I(I,\phi )=\gamma \log (J(\phi ))+I\text{ 
}(\func{mod}\text{ }\pi ) 
\end{array}
,\right. 
$$
that is also induced by an element of $TQ.$

And from \cite{ragazzo}, $\exists $ $\alpha _{crit}(\gamma )$ such that for $%
\alpha >\alpha _{crit}(\gamma ),$ $\overline{F}$ does not have R.I.C's. In
this case, we can apply the same criteria explained for the standard map.
But as there are 2 parameters, we do not obtain a critical value, we obtain
a critical set in the $(\gamma ,\alpha )$ plane. Another important
application of this theory is to obtain properties about the structure of
the unstable set of the above mentioned homoclinic loop (to the
saddle-center equilibrium), when the former is unstable. The periodic orbits
given by theorem (\ref{PtosPerQuoc}) were analyzed in \cite{zanrag1} and it
was proved that for every vertical rotation number $\frac mn>0$, there is an
open set in the parameter space with a $\frac mn$-periodic orbit which is
topologically a sink. In particular, it can be proved that, for a fixed
value of $\gamma >0,$ given an $\epsilon >\alpha _{cr}(\gamma )>1,$ where $%
\alpha _{cr}(\gamma )$ is analogous to the constant $k_{cr}$ defined for the
standard map, there is a number $\frac mn>0$ and an open interval $I_{\frac
mn}\subset (\alpha _{cr}(\gamma ),\epsilon ),$ such that for $\alpha \in
I_{\frac mn},$ $\overline{F}$ has a vertical periodic orbit with $\rho
_V=\frac mn$ which is also a topological sink. So we can say that one of the
mechanisms that cause the lost of stability of the homoclinic loop is the
creation of periodic sinks for $\overline{F}$. And in \cite{zanrag2} it was
proved that the existence of a topological sink for $\overline{F}$ implies
many interesting properties on the topology of the set of orbits that have
the saddle-center loop as their $\alpha $-limit set (a set analogous to the
unstable manifold of a hyperbolic periodic orbit). More precisely, in this
case, given an arbitrary neighborhood of the original homoclinic loop, a set
of positive measure contained in this neighborhood escapes from it following
\ (or clustering around) a finite set of orbits that in a certain sense,
correspond to the topological sinks for $\overline{F}.$ In a forthcoming
paper, we will analyze the following function:

$$
\rho _V^{\max }(\gamma ,\alpha )=\stackunder{P\in {\rm T^2}}{\sup }\rho
_V(P)=\stackunder{P\in {\rm T^2}}{\sup }\left[ \stackunder{n\rightarrow
\infty }{\lim }\frac{p_2\circ F^n(P)-p_2(P)}n\right] ,\text{ } 
$$
where the supremum is taken over all $P\in {\rm T^2}$ such that $\rho _V(P)$
exists. Using a method developed in \cite{zanrag1} and results from \cite
{misiuzie}, we plan to prove the density of periodic sinks in the subset of
the parameter space $(\gamma ,\alpha )$ where $\overline{F}$ does not have
R.I.C's.

\vskip0.3truecm

3) Given a $C^2$ circle diffeomorphism $f:S^1\rightarrow S^1$ $(f(\phi
+1)=f(\phi )+1)$ we can define the following generating function: 
$$
h_f(\phi ,\phi ^{\prime })=\exp (\phi ^{\prime }-f(\phi ))\text{ } 
$$

As $h_f(\phi +1,\phi ^{\prime }+1)=h_f(\phi ,\phi ^{\prime })$ the
associated twist map $\widehat{T}_f:S^1\times ]0,+\infty [\hookleftarrow $
is exact: 
$$
\widehat{T}_f:\left\{ 
\begin{array}{l}
\phi ^{\prime }=f\;(\phi )+\log (I^{\prime }) 
\text{ }(\func{mod}\text{ }1) \\ I^{\prime }=\frac 1{f^{\prime }(\phi )}I 
\end{array}
\right. 
$$

By the same coordinate change applied to $\widehat{F}$ 
$$
\left\{ 
\begin{array}{l}
\stackrel{\symbol{126}}{\phi }=\phi \\ \stackrel{\symbol{126}}{I}=\log (I) 
\end{array}
\right. 
$$

we can write $\widehat{T}_{f}$ in the following way: 
$$
\left\{ 
\begin{array}{l}
\phi ^{\prime }=f\;(\phi )+I^{\prime } 
\text{ }(\func{mod}\text{ }1) \\ I^{\prime }=\log (\frac{1}{f^{\prime }(\phi
)})+I 
\end{array}
\right. 
$$

As above, in these coordinates $\widehat{T}_f$ is $\mu $-exact for the
following measure: 
$$
\mu (A)=\stackunder{A}{\int }e^Id\phi dI 
$$

And $\widehat{T}_f$ induces a torus map 
$$
\overline{T_f}:\left\{ 
\begin{array}{l}
\phi ^{\prime }=f\;(\phi )+I^{\prime } 
\text{ }(mod\text{ }1) \\ I^{\prime }=\log (\frac 1{f^{\prime }(\phi )})+I 
\text{ }(mod\text{ }1) 
\end{array}
\right. 
$$
such that our results apply.

\vskip0.2truecm

{\it Acknowledgements: }I am very grateful to C. Grotta Ragazzo for
listening to oral expositions of these results, reading the first
manuscripts and for many discussions, comments and all his support, to J.
Mather for all his support and to the referees for a very careful reading of
the paper, for the suggestion of a new proof of theorem (\ref{PtosPerQuoc}),
for comments on how to obtain 2 periodic orbits in theorem (\ref
{ImplicPeriod}) in the general case and for all their other remarks that
improved the text.


\end{document}